\documentclass[12pt]{article}
\usepackage[dvipsnames]{xcolor}
\usepackage[obeyspaces,hyphens,spaces]{url}

\usepackage[utf8,latin1]{inputenc}
\usepackage[T1]{fontenc}
\usepackage{tikz}
\usepackage{pgfplots}
\pgfplotsset{compat=1.9}
\usepackage[percent]{overpic}
\usepackage{xfrac}
\usepackage{amsmath}
\usepackage{float}
\usepackage{soul}
\usepackage{mathpazo}
\usepackage{empheq}
\usepackage{geometry}
\usepackage[shortlabels,inline]{enumitem}
\usepackage{booktabs}

\usepackage[doc]{optional}

\usepackage{float}
\usepackage{soul}
\usepackage{graphicx}
\definecolor{labelkey}{rgb}{0,0.08,0.45}
\definecolor{refkey}{rgb}{0,0.6,0.0}
\definecolor{Brown}{rgb}{0.45,0.0,0.05}
\definecolor{lime}{rgb}{0.00,0.8,0.0}
\definecolor{lblue}{rgb}{0.5,0.5,0.99}

\usepackage{amsmath}
\usepackage{amssymb}
\usepackage{stmaryrd}
\usepackage{pgfplots}
\pgfplotsset{compat=1.6}
\definecolor{dgreen}{rgb}{0.00,0.49,0.00}
\definecolor{dblue}{rgb}{0,0.08,0.55}
\colorlet{minhblue}{dblue}
\colorlet{mygreen}{dgreen}
\usepackage{hyperref}
\hypersetup{ 
    linktocpage=true,
    colorlinks= true, 
    urlcolor = Red, 
    citecolor = mygreen,
    linkcolor = minhblue
}

\newcommand{\thmtit}[1]{{\bf{#1}}}
\newcommand*{\tran}{^{\mkern-1.5mu\mathsf{T}}}
\providecommand{\siff}{\Leftrightarrow}

\newcommand{\nnn}{\ensuremath{{n\in{\mathbb N}}}}

\newcommand{\menge}[2]{\big\{{#1}~\big |~{#2}\big\}}

\newcommand{\fenv}[1]%
{\ensuremath{\,\overrightarrow{\operatorname{env}}_{#1}}}
\newcommand{\benv}[1]%
{\ensuremath{\,\overleftarrow{\operatorname{env}}_{#1}}}

\newcommand{\RR}{\ensuremath{\mathbb R}}

\newcommand{\NN}{\ensuremath{\mathbb N}}

\newcommand{\dom}{\ensuremath{\operatorname{dom}}}
\newcommand{\argmin}{\ensuremath{\operatorname{argmin}}}

\newcommand{\zer}{\ensuremath{\operatorname{zer}}}

\newcommand{\Id}{\ensuremath{\operatorname{Id}}}

\newcommand{\bx}{\ensuremath{\mathbf{x}}}
\newcommand{\by}{\ensuremath{\mathbf{y}}}

\newcommand{\bA}{\ensuremath{{\mathbf{A}}}}
\newcommand{\bB}{\ensuremath{{\mathbf{B}}}}

{\begin{list}{}{%
\settowidth{\labelwidth}{\textrm{#1~}}%
\setlength{\leftmargin}{\labelwidth+\labelsep}}}%requires macro calc.sty
{\end{list}}
\usepackage{amsthm}
\usepackage[capitalize, nameinlink]{cleveref}
\crefname{equation}{}{equations}
\crefname{chapter}{Appendix}{chapters}
\crefname{item}{}{items}
\newtheorem{theorem}{Theorem}[section]
\newtheorem{lemma}[theorem]{Lemma}

\newtheorem{corollary}[theorem]{Corollary}

\newtheorem{proposition}[theorem]{Proposition}

\newtheorem{example}[theorem]{Example}
\newtheorem{ex}[theorem]{Example}

\newtheorem{remark}[theorem]{Remark}

\providecommand{\norm}[1]{\lVert#1\rVert}
\providecommand{\normsq}[1]{\lVert#1\rVert^2}

\providecommand{\stb}[1]{\left\{#1\right\}}
\providecommand{\innp}[1]{\langle#1\rangle}

\providecommand{\RA}{\Rightarrow}

\providecommand{\grad}{\nabla}

\providecommand{\RR}{\mathbb{R}}

\providecommand{\dom}{\operatorname{dom}}

\newcommand{\fix}{\ensuremath{\operatorname{Fix}}}

\providecommand{\gra}{\operatorname{gra}}
\providecommand{\Id}{\operatorname{{ Id}}}

\providecommand{\fady}{\varnothing}

\providecommand{\argmin}{\mathrm{arg}\!\min}

\providecommand{\rras}{\rightrightarrows}

\providecommand{\NN}{\mathbb{N}}

\providecommand{\fix}{\operatorname{Fix}}

\providecommand{\Id}{\operatorname{Id}}

\providecommand{\pt}{{\partial}}

\providecommand{\zer}{\operatorname{zer}}
\providecommand{\R}{{ R}}

\providecommand{\fady}{\varnothing}

\providecommand{\DR}{\text{\scriptsize DR}}

\providecommand{\RR}{\mathbb{R}}
\providecommand{\NN}{\mathbb{N}}

\definecolor{myblue}{rgb}{.8, .8, 1}

\allowdisplaybreaks % or locally if problems {\allowdisplaybreaks
\begin{document}
%-------------------------------------------------------------------------

%\tikzstyle{decision} = [diamond, draw, fill=blue!50]
%\tikzstyle{line} = [draw, -stealth, thick]
%\tikzstyle{elli}=[draw, ellipse, fill=red!50,minimum height=8mm, text width=5em, text centered]
%\tikzstyle{block} = [draw, rectangle, fill=blue!50, text width=8em, text centered, minimum height=15mm, node distance=10em]
%

\title{Douglas--Rachford splitting
for a Lipschitz continuous and a strongly monotone operator}
\author{Walaa M.\ Moursi\thanks{Department of 
Electrical Engineering,
Stanford 
University, %Packard Building, Room 222,
% 350 Serra Mall, 
 Stanford, CA 94305, USA
%and 
%Mansoura University, Faculty of Science, Mathematics Department, 
%Mansoura 35516, Egypt. 
E-mail:
\texttt{wmoursi@stanford.edu}.}
~and Lieven Vandenberghe\thanks{Department of 
Electrical and Computer Engineering, 
University of California Los Angeles, Los Angeles, CA 90095, USA .
E-mail: \texttt{vandenbe@ucla.edu}.}}

%\date{November 16, 2017}

\maketitle

\begin{abstract}
\noindent
The Douglas--Rachford method
is a popular splitting technique 
for finding a zero of the sum of
two subdifferential operators of proper 
closed convex functions;
more generally
two maximally monotone operators.
Recent results concerned with linear rates of convergence of the
method require additional properties
 of the underlying monotone operators,
such as strong monotonicity and cocoercivity. 
In this paper, we study the case when 
one operator is Lipschitz continuous
but not necessarily a subdifferential operator 
and the other operator is strongly monotone.
This situation arises in optimization methods 
which involve primal-dual approaches.
We provide new linear convergence results in this setting.

\end{abstract}
{\small
\noindent
{\bfseries 2010 Mathematics Subject Classification:}
{
Primary 
%47H05, %Monotone operators and generalizations
%47H09, %Contraction-type mappings, nonexpansive mappings, $A$-proper mappings,
%49M27; %Decomposition Methods
%Secondary 
%49M29, %Methods involving duality etc.
%49N15, %Duality theory
%90C25. 
47H05, %Monotone operators and generalizations
47H09, %Contraction-type mappings, nonexpansive mappings, $A$-proper mappings,
49M27; %Decomposition Methods,
90C25. 
Secondary 
49M29, %Methods involving duality etc.
49N15. %Duality theory

}

\noindent {\bfseries Keywords:}
%Attouch--Th\'era duality, 
%Alternating Direction Method of Multipliers (ADMM),
cocoercive operator,
Douglas--Rachford algorithm,
%Dykstra method,
%Fenchel--Rockafellar Duality,
%Method of Alternating Projections (MAP),
%Peaceman--Rachford algorithm,
%inconsistent case,
linear convergence,
Lipschitz continuous mapping, 
maximally monotone operator,
nonexpansive mapping,
rate of convergence,
skew symmetric operator,
splitting methods,
strongly convex function,
strongly monotone operator.
%paramonotone operator,
%sum problem,
%weak convergence.
}

%\section{Introduction and motivation}

%All acknowledgements should be placed in the back of the paper after Conclusions..

\section{Introduction}

The Douglas--Rachford splitting algorithm, 
introduced by Lions and
Mercier~\cite{L-M79}, is a fundamental algorithm for solving 
monotone inclusion problems that involves finding a zero
of the sum of two maximally monotone operators $A$ and $B$.
%\begin{equation}
% \label{eq:P:sum}
% \text{Find $x\in X$ such that $0\in Ax+Bx$,}
%\end{equation}
%where  $X$ is a real Hilbert space, and 
%$A\colon X\rras X$ and $B\colon X\rras X$ are maximally 
%monotone operators.   
(See Section~\ref{sec:RA:prop} for a review of these
and other definitions used in the paper.)
%The algorithm is based on the iteration
%\begin{equation}
% u_{n+1} = T_{\DR}u_n
% = \tfrac{1}{2} (\Id + R_BR_A)u_n,
%\end{equation}
%starting at arbitrary $u_0\in X$, where $\Id$ is the identity operator, 
%and $R_A$ and $R_B$ are the 
%\emph{reflected resolvents} of $A$ and $B$
%(defined as $R_A = 2J_A - I$ and $R_B = 2J_B-I$  where  the operators
%$J_A = (\Id + A)^{-1}$ and $J_B = (\Id + B)^{-1}$ are the 
%\emph{resolvents}).
Monotone inclusions can be used to formulate primal, dual, and
primal-dual optimality conditions of convex optimization problems,
equilibrium conditions in convex-concave games,
monotone variational inequalities, 
and monotone complementarity problems.
The Douglas--Rachford algorithm is useful for all these applications, 
provided that the operator in the inclusion 
problem can be written as a sum 
of two operators, as in \cref{eq:P:sum}, with resolvents that are easily 
computed.
This is often the case in large scale applications; see, e.g., 
\cite{BC2017,Borwein50,Brezis,BurIus,Comb96,Simons1,Simons2,Rock98,Zeidler2a,Zeidler2b}, and the references therein.  
The Douglas--Rachford method can also be used to derive other important 
splitting methods, such as the 
Alternating Direction Method of Multiplier or ADMM \cite{Gabay83,EckBer},
Spingarn's method of partial inverses \cite{EckBer},
the primal-dual hybrid gradient method \cite{OV17}, 
and linearized ADMM \cite{OV17}.

Under additional assumptions on the operators 
$A$ and $B$, linear rates of convergence are possible.
In their seminal work \cite{L-M79}, Lions and Mercier 
proved linear convergence of the Douglas--Rachford
iteration when one operator is strongly monotone and 
cocoercive.
%(i.e.,  both the operator and its inverse are strongly monotone).
Recent works concerned with linear
rates of convergence of Douglas--Rachford method 
include \cite{JAT2014,BM15,Phan16}
for linear rates in convex feasibility settings, 
\cite{Gis17,GB17,L-M79} for linear rates 
under strong convexity assumptions,
\cite{DZ14} for linear rates in basis pursuit setting, 
and \cite{ACL16,HL13,HLN14} for local linear rates in more general 
settings.
In the recent work \cite{Gis17},
Giselsson studied and proved tight linear rates of convergence 
of Douglas--Rachford in the following three cases:
\begin{enumerate*}[(a)]
\item
\label{case:a}
$A$ is strongly monotone and $B$
is cocoercive
(see \cite[Theorem~5.6]{Gis17}),
\item
\label{case:b}
$A$ is strongly monotone and Lipschitz continuous
(see \cite[Theorem~6.5]{Gis17}),
and 
\item
\label{case:c}
$A$ is strongly monotone and cocoercive
(see \cite[Theorem~7.4]{Gis17}).
Giselsson's results are independent of the order 
of $A$ and $B$, and therefore actually cover \emph{six} cases.
\end{enumerate*}

The main contribution of this paper is to supplement Giselsson's results 
with a linear convergence result for the case when
$A$ is Lipschitz continuous and $B$ is strongly monotone.
Unlike the results in~\cite{Gis17}, our linear convergence result
is not symmetric in $A$ and $B$, and does not apply to the case where 
$A$ is strongly monotone and $B$ is Lipschitz continuous,
except in the important special case when $B$ is a linear mapping.
When $A$ is the subdifferential of a convex function, 
Lipschitz continuity and cocoercivity are equivalent properties.  
However, for general monotone operators, 
Lipschitz continuity is a much weaker condition than cocoercivity,
so the case studied in this paper is an important extension of 
\cite[Theorem~5.6]{Gis17}.

As an application, we discuss the Douglas--Rachford splitting method
applied to the primal-dual optimality conditions of a convex
problem, formulated as an inclusion problem 
\cref{eq:P:sum} in which one of
the operators is a skew-symmetric linear mapping and  \emph{not} 
a subdifferential, see, e.g.,  \cite{BH13,B-AC11,CP11,Cond13,OV14}.

This paper is organized as follows.
Section~\ref{sec:RA:prop} presents a collection 
of useful properties of reflected resolvents
of monotone operators under 
additional assumptions on the operator.
Section~\ref{sec:Motivation} provides a high level overview of 
relevant linear convergence results.
Our main results appear in
Section~\ref{sec:Main} where 
we prove linear convergence of 
Douglas--Rachford iteration when 
applied to find a zero of the sum of 
maximally monotone operators $A $ and $B$ when $A$ is
Lipschitz continuous 
and
$B$ is strongly monotone.
Finally, in Section~\ref{sec:Applications}
we present an application of our results to the 
primal dual Douglas--Rachford method.

%\subsection*{Notation}

\section{Contraction properties of reflected resolvents}
\label{sec:RA:prop}
Throughout the paper, $X$ is a real Hilbert space,
with inner product $\innp{\cdot,\cdot}$
and induced norm $\norm{\cdot}$.
We use  the notation $A\colon X\rras X$ to indicate that
$A$ is a set-valued operator on $X$.
The \emph{domain} of $A$ is 
$\dom (A)=\menge{x\in X}{Ax\neq \varnothing}$
and the \emph{graph} of $A$
is $\gra (A)=\menge{(x,u)\in X\times X}{u\in Ax}$.
We use the notation $A\colon X\to X$ to indicate 
that A is a single-valued operator on $X$ and $\dom A=X$. 
The inverse of $A$, denoted by $A^{-1}$, is the operator
with graph 
$\gra (A^{-1}) = \{(u,x) \mid (x,u) \in \gra(A)\}$.
Let $C$ be a convex closed nonempty subset  of $X$.
We use $N_C$ to denote
the \emph{normal cone} of $C$ defined as
$N_C(x)\coloneqq\menge{u\in X }{\innp{u,y-x}\le 0 
\text{  for all  }y\in C}$, if $x\in C$; and $N_C(x)\coloneqq\fady$, otherwise;
and $P_C$ to denote that 
orthogonal projection onto $C$
(this is also known as the closest point mapping)
defined at every $x\in X$ by 
$P_C(x)\coloneqq \argmin_{c\in C}\norm{x-c}$. 

An operator $A$ on $X$ is $\beta$-\emph{Lipschitz continuous} if 
it is single-valued on $\dom(A)$ and 
\begin{equation}
  \norm{Ax-Ay}\le \beta \norm{x-y} \quad \forall x, y \in \dom(A).
\end{equation}
A $1$-Lipschitz continuous operator is called \emph{nonexpansive}.
An operator $A$ is $\alpha$-\emph{averaged}, with 
$\alpha\in\left[0,1\right[$, if it can be expressed as
$A=(1-\alpha)\Id+\alpha N$ where $\Id$ is the identity operator,
 $N$ is nonexpansive.
An operator $A\colon X\rras X$ is \emph{monotone} if 
\begin{equation}
 \innp{x - y, u-v} \geq 0 \quad \forall (x,u), (y,v) \in \gra(A).
\end{equation}
A monotone operator $A$ is \emph{maximally monotone} if its graph 
admits no proper extension that preserves the monotonicity of $A$.
An operator $A$ is $\mu$-\emph{strongly monotone}, with $\mu>0$,
if
\begin{equation}
 \innp{x - y, u-v} \ge \mu \norm{x-y}^2 \quad
 \forall (x,u), (y,v) \in \gra(A).
\end{equation}
Equivalently, $A-\mu\Id$ is monotone.
The operator $A$ is $(1/\beta)$-\emph{cocoercive}, with $\beta > 0$,
if its inverse $A^{-1}$  is
$(1/\beta)$-strongly monotone, i.e.,
\begin{equation}
\label{eq:Lip:coco}
 \innp{x - y, u-v} \ge \tfrac{1}{\beta} \norm{u-v}^2 \quad
 \forall (x,u), (y,v) \in \gra(A).
\end{equation}
Note that this implies that $A$ is single-valued on its domain,
and (by the Cauchy-Schwarz inequality) that $A$ is 
$\beta$-Lipschitz continuous.
A $1$-cocoercive operator is also called \emph{firmly nonexpansive}.
We note that all these properties are defined as 
quadratic inequalities on the graph of the operator.
Table~\ref{tab:1} summarizes the definitions.
Each of the four properties in the table is defined as  
\begin{equation}
\label{eq:quad:ineq}
L_{11} \norm{x-y}^2 + 2 L_{12} \innp{x-y, u-v}
 + L_{22} \norm{u-v}^2\ge 0 \quad \forall (x,u), (y,v) \in \gra(A),
\end{equation}
where $L$ is the $2\times 2$ matrix shown on row~2 of the table.
%We assume $\mu>0$, $\beta>0$, and $\alpha\in \left]0,1\right[$.

\begin{table}
\centering \small
\begin{tabular}{
@{}c@{\hskip 1em}c@{\hskip 2em}c@{\hskip 2em}c@{\hskip 2em}c@{}} \toprule
& $\mu$-strong monotonicity
& $\beta$-Lipschitz continuity
& $(1/\beta)$-cocoercivity
& $\alpha$-averagedness\\ \cmidrule{2-5} 
\addlinespace[\belowrulesep]
$L$ & 
$\begin{bmatrix} -2\mu&1\\ 1&0 \end{bmatrix}$&
$\begin{bmatrix} \beta^2&0\\ 0&-1 \end{bmatrix}$&
$\begin{bmatrix} 0&\beta\\ \beta&-2 \end{bmatrix}$&
$\begin{bmatrix} 2\alpha-1&1-\alpha\\ 1-\alpha&-1 \end{bmatrix}$\\
\\
$M$ & 
$\begin{bmatrix} 0&1\\ 1&-2\mu-2 \end{bmatrix}$&
$\begin{bmatrix} -1&1\\ 1&\beta^2-1 \end{bmatrix}$&
$\begin{bmatrix} -2&\beta+2\\ \beta+2&-2\beta-2 \end{bmatrix}$&
$\begin{bmatrix} -1&2-\alpha\\ 2-\alpha&4\alpha-4 \end{bmatrix}$\\
\\
$N$&
$2\begin{bmatrix} 1-\mu&-\mu\\ -\mu&-1-\mu \end{bmatrix}$& 
$\begin{bmatrix} \beta^2-1&\beta^2+1\\ \beta^2+1&\beta^2-1
\end{bmatrix}$&
$\displaystyle
\frac{2}{\beta} \begin{bmatrix} \beta-1&1\\ 1&-\beta-1 \end{bmatrix}$&
$2\begin{bmatrix} 0&\alpha\\ \alpha&-2(1-\alpha) \end{bmatrix}$
\\
\addlinespace[\belowrulesep]
\bottomrule
\end{tabular}
\caption{Each of the four operator properties is defined 
as \cref{eq:quad:ineq} for the matrix $L$ shown in the table.
They can be defined equivalently as properties of the resolvent,
given by \cref{eq:quad:ineq:JA} with the matrix $M$ shown in the table,
and as properties of the reflected resolvent, 
given by \cref{eq:quad:ineq:RA} for the matrix $N$ shown in the table.  
Taking $\mu=0$ in the first column also gives three equivalent 
definitions of monotonicity.}
\label{tab:1}
\end{table}

The \emph{resolvent} of an operator $A$ is the mapping 
$J_A = (\Id +A)^{-1}$.  The \emph{reflected resolvent} is the mapping
$R_A = 2J_A -\Id$.
The graphs of the resolvent and reflected resolvent 
of an operator $A$
are related to the graph of $A$ by invertible linear transformations:
\begin{align}
\label{eq:A:JA}
\gra A
& =  \{ (u, x-u) \mid (x,u)\in \gra(J_A)\}  \\
\label{eq:A:RA}
& =  \left\{ \tfrac{1}{2} (x+u, x-u) \mid (x,u)\in \gra(R_A)\right\}.
\end{align}
Hence, if we define
\begin{equation}
\label{eq:LMN}
M=
\begin{bmatrix}
0&1\\
1&-1
\end{bmatrix}
L
\begin{bmatrix}
0&1\\
1&-1
\end{bmatrix},
\quad
N=
\begin{bmatrix}
1&1\\
1&-1
\end{bmatrix}
L
\begin{bmatrix}
1&1\\
1&-1
\end{bmatrix},
\end{equation}
then the property~\cref{eq:quad:ineq} is equivalent
to
\begin{equation}
\label{eq:quad:ineq:JA}
M_{11}\norm{x-y}^2
+2M_{12}\innp{x-y,u-v}
+M_{22}\norm{u-v}^2
\ge 0\quad  \forall (x,u), (y,v) \in \gra(J_A),
\end{equation} 
and also to
\begin{equation}
\label{eq:quad:ineq:RA}
N_{11}\norm{x-y}^2
+2N_{12}\innp{x-y,u-v}
+N_{22}\norm{u-v}^2
\ge 0\quad  \forall (x,u), (y,v) \in \gra(R_A).
\end{equation} 
For each property in~Table~\ref{tab:1}, we therefore have
equivalent definitions of the 
form~\cref{eq:quad:ineq:JA}~and~\cref{eq:quad:ineq:RA}.  The matrices 
$M$ and $N$ are shown in the third and fourth rows of the table,
respectively.

In \cref{prop:corres:A:RA} below we collect some
useful properties of the reflected resolvent $\R_A$.
We point out that
items \ref{prop:corres:A:RA:ii} \& \ref{prop:corres:A:RA:iii}
below provide short proofs for 
Theorems~7.2~and~6.3 in \cite{Gis17}. 
\begin{proposition}
\label{prop:corres:A:RA}
The following statements hold for any operator $A$.
\begin{enumerate}
\item
\label{prop:corres:A:RA:0}
Suppose $\mu> 0$ and $\beta > 0$.
If $A$ is $\mu$-strongly monotone
and $\beta$-Lipschitz continuous, then
$A$ is $(\mu/\beta^2)$-cocoercive.
\item
\label{prop:corres:A:RA:v}
$A$ is monotone and nonexpansive 
if and only if $J_A$ is $\tfrac{1}{2}$-strongly monotone,
and if and only if $R_A$ is monotone.
\item 
\label{prop:corres:A:RA:vii}
Suppose $\mu > 0$.
$A$ is $\mu$-strongly monotone 
if and only if 
$J_A$ is $(1+\mu)$-cocoercive, and if and only if 
$-R_A$ is $(1+\mu)^{-1}$-averaged.
\item
\label{prop:corres:A:RA:ii}
Suppose $\beta \geq  \mu > 0$.
If $A$ is $\mu$-strongly monotone and $(1/\beta)$-cocoercive,
then $R_A$ is $\kappa$-Lipschitz continuous with
\begin{equation}
\kappa=\left(\frac{1-2\mu+\mu\beta}{1+2\mu+\mu\beta}\right)^{1/2}.
\end{equation}

\item
\label{prop:corres:A:RA:iii}
Suppose $\beta \geq \mu > 0$.
If $A$ is $\mu$-strongly monotone and $\beta$-Lipschitz
continuous, then $R_A$ is $\kappa$-Lipschitz continuous with
\begin{equation}
\label{e-kappa-strongly-monotone-lipschitz}
\kappa=\left(\frac{1-2\mu+\beta^2}{1+2\mu+\beta^2}\right)^{1/2}.
\end{equation}

\item
\label{prop:corres:A:RA:i}
Suppose $0 < \mu < 1$ and $0 < \alpha < 1$.
If $A$ is $\mu$-strongly monotone and $\alpha$-averaged,
then $R_A$ is $\kappa$-Lipschitz continuous with
\begin{equation}
\kappa=\left(\frac{\alpha(1-\mu)}{\alpha(1-\mu)+2\mu}\right)^{1/2}.
\end{equation}

\end{enumerate}
\end{proposition}

{\it Proof}
\ref{prop:corres:A:RA:0}:
This follows from Table~\ref{tab:1} and
\begin{equation}
\begin{bmatrix}
0&\beta^2/\mu\\
\beta^2/\mu&-2
\end{bmatrix}
=
\frac{\beta^2}{\mu} \begin{bmatrix}
-2\mu&1\\
1&0
\end{bmatrix}
+
2\begin{bmatrix}
\beta^2&0\\
0&-1
\end{bmatrix}.
\end{equation}
\ref{prop:corres:A:RA:v}:
Use Table~\ref{tab:1}
for the case when $A$
is $1$-Lipschitz continuous.
\ref{prop:corres:A:RA:vii}:
This is clear from Table~\ref{tab:1}. 
\ref{prop:corres:A:RA:ii}:
This follows from Table~\ref{tab:1} and
\begin{equation}
\begin{bmatrix}
1-2\mu+\mu\beta&0\\
0&-(1+2\mu+\mu\beta)
\end{bmatrix}
=
\begin{bmatrix}
1-\mu&-\mu\\
-\mu&-1-\mu
\end{bmatrix}
+\mu
\begin{bmatrix}
\beta-1&1\\
1&-\beta-1
\end{bmatrix}.
\end{equation}
\ref{prop:corres:A:RA:iii}:
This follows from Table~\ref{tab:1} and
\begin{equation} \label{e-kappa-strongly-monotone-lipschitz-pf}
\begin{bmatrix}
1-2\mu+\beta^2&0\\
0&-(1+2\mu+\beta^2)
\end{bmatrix}
=
(\beta^2+1)
\begin{bmatrix}
1-\mu&-\mu\\
-\mu&-1-\mu
\end{bmatrix}
+\mu
\begin{bmatrix}
\beta^2-1&\beta^2+1\\
\beta^2+1&\beta^2-1
\end{bmatrix}.
\end{equation}
Alternatively, combine \ref{prop:corres:A:RA:0}
 and \ref{prop:corres:A:RA:ii} 
 applied with $\beta$ replaced
 by $\beta^2/\mu$ to learn that 
 $R_A$
is $\kappa$-Lipschitz continuous with
\begin{equation}
\kappa=\left(\frac{1-2\mu+\mu(\beta^2/\mu)}{1+2\mu+\mu(\beta^2/\mu)}\right)^{1/2}
=\left(\frac{1-2\mu+\beta^2}{1+2\mu+\beta^2}\right)^{1/2}.
\end{equation}
\ref{prop:corres:A:RA:i}:
This follows from Table~\ref{tab:1} and the identity
\begin{equation}
\begin{bmatrix}
\alpha(1-\mu)&0\\
0&-(\alpha(1-\mu)+2\mu)
\end{bmatrix}
=
\alpha\begin{bmatrix}
1-\mu&-\mu\\
-\mu&-1-\mu
\end{bmatrix}
+\mu
\begin{bmatrix}
0&\alpha\\
\alpha&-2(1-\alpha)
\end{bmatrix}.
\end{equation}
% \ref{prop:corres:A:RA:iv}:
% Apply either \ref{prop:corres:A:RA:ii} with $\beta=1$
%or \ref{prop:corres:A:RA:i} with $\alpha=1/2$.
\qed

\begin{remark}
The contraction factors of the reflected resolvents 
are important in the linear convergence proofs in \cite{Gis17}.
\cref{prop:corres:A:RA}\ref{prop:corres:A:RA:iii} gives  
the contraction factor of the
reflected resolvent of a strongly monotone and Lipschitz
continuous operator.   As indicated in the proof, this 
result can be derived in two ways.  In the second approach, we use
\cref{prop:corres:A:RA}\ref{prop:corres:A:RA:0} to derive
the contraction factor from the result for strongly monotone and
cocoercive operators. 
\end{remark}
We conclude this section with 
the following lemma.
\begin{lemma}
\label{lem:RA:Lips:beta}
Suppose that $A\colon X\to X$ is monotone and 
$\beta$-Lipschitz continuous with $\beta>0$. 
Let $(x,y)\in X\times X$.
Then the following hold:
%\begin{equation}\end{equation}
\begin{enumerate}
\item
\label{lem:RA:Lips:beta:i}
$\norm{x-y}\le(1+\beta) \norm{J_A x-J_Ay}$. 
%is $\tfrac{1}{(\beta+1)^2}$-strongly monotone.

\item
\label{lem:RA:Lips:beta:ii}
$\Id-J_A$ is a Banach contraction with constant 
$\tfrac{\beta}{\sqrt{1+\beta^2}}$. 
\item
\label{lem:RA:Lips:beta:iii}
$J_A$ is $\left(\tfrac{1}{2(1+\beta)^2}
+\tfrac{1}{2(1+\beta^2)}\right)$-strongly
monotone.

\item
\label{lem:RA:Lips:beta:iv}
$\innp{x-y,R_Ax-R_Ay}
\ge- \lambda
\norm{x-y}^2$ where\footnote{This property also means that 
$R_A$ is \emph{hypomonotone}, see \cite[Example~12.28]{Rock98}. }
$\lambda=\left(1-\tfrac{1}{(1+\beta)^2}-\tfrac{1}{1+\beta^2}\right)\in 
\left]-1,1\right[$.

\end{enumerate}
\end{lemma}
{\it Proof}
\ref{lem:RA:Lips:beta:i}:
This follows from 
entry $(2,2)$ of Table~\ref{tab:1}
and 
\begin{equation}
\begin{bmatrix}
-1&0\\
0&(1+\beta)^2
\end{bmatrix}
\succeq
\tfrac{\beta+1}{\beta}
\begin{bmatrix}
-1&1\\
1&\beta^2-1
\end{bmatrix}.
\end{equation}

\ref{lem:RA:Lips:beta:ii}:
Let $P$ be a $2\times 2$ matrix 
satisfying
\begin{equation}
\label{eq:quad:ineq:JA:inv}
P_{11}\norm{x-y}^2
+2P_{12}\innp{x-y,u-v}
+P_{22}\norm{u-v}^2
\ge 0\quad  \forall (x,u), (y,v) \in \gra(J_{A^{-1}}).
\end{equation} 
On the one hand, it follows from
\cref{eq:quad:ineq}
that each of the four properties in the 
first row of Table~\ref{tab:1}   
correspond to 
\begin{equation}
\label{eq:quad:ineq:invA}
L_{22} \norm{x-y}^2 + 2 L_{12} \innp{x-y, u-v}
 + L_{11} \norm{u-v}^2\ge 0 \quad \forall (x,u), (y,v) \in \gra(A^{-1}).
\end{equation}
Therefore, using 
the $(1,2)$
entry of Table~\ref{tab:1},
the first equation in 
\cref{eq:LMN} (applied to $A^{-1}$)
 and 
\cref{eq:quad:ineq:invA}
we learn that 
the $\beta$-Lipschitz continuity of $A$
corresponds to the matrix
\begin{equation}
\label{eq:P:Lips}
P=
\begin{bmatrix}
0&1\\
1&-1
\end{bmatrix}
\begin{bmatrix}
-1&0\\
0&\beta^2
\end{bmatrix}
\begin{bmatrix}
0&1\\
1&-1
\end{bmatrix}
=
\begin{bmatrix}
\beta^2&-\beta^2\\
-\beta^2&\beta^2-1
\end{bmatrix}.
\end{equation}
Similarly, we learn from 
the $(2,1)$
entry of Table~\ref{tab:1} (applied with $\mu=0$),
that the monotonicity of $A$ (equivalently, the monotonicity
of $A^{-1}$) corresponds to 
the matrix
\begin{equation}
\label{eq:P:mono}
P=
\begin{bmatrix}
0&1\\
1&-2
\end{bmatrix}.
\end{equation}

The conclusion then follows 
from \cref{eq:P:Lips}
and
\cref{eq:P:mono}
in view of 
\cref{eq:quad:ineq:JA:inv}
by noting that
\begin{equation}
\begin{bmatrix}
\beta^2&0\\
0&-(1+\beta^2)
\end{bmatrix}
=
\begin{bmatrix}
\beta^2&-\beta^2\\
-\beta^2&\beta^2-1
\end{bmatrix}
+
\beta^2
\begin{bmatrix}
0&1\\
1&-2
\end{bmatrix}.
\end{equation}

\ref{lem:RA:Lips:beta:iii}:
This follows from the entries $(2,1)$
(applied with $\mu=0$)
and $(2,2)$ in Table~\ref{tab:1},
\ref{lem:RA:Lips:beta:i}
 and 
 \begin{align}
& \begin{bmatrix}
-1/(1+\beta^2)-1/(1+\beta)^2&1\\
1&0
\end{bmatrix}
\nonumber
\\
&\quad=
\tfrac{1}{\beta^2+1}
\begin{bmatrix}
-1&1\\
1&\beta^2-1
\end{bmatrix}
+
\tfrac{\beta^2}{1+\beta^2}
\begin{bmatrix}
0&1\\1&-1
\end{bmatrix}
+
\tfrac{1}{(1+\beta^2)}
\begin{bmatrix}
-1&0\\
0&(1+\beta)^2
\end{bmatrix},
 \end{align}
 in view of 
\cref{eq:quad:ineq:JA:inv}.

\ref{lem:RA:Lips:beta:iv}:
One can readily verify that
\begin{equation}
\label{eq:graJA:RA}
\gra(J_A)=\tfrac{1}{2}
\begin{bmatrix}
2&0\\1&1
\end{bmatrix}
\gra (R_A).
\end{equation}
Now the conclusion follows from 
\cref{eq:graJA:RA}
and 
\ref{lem:RA:Lips:beta:iii}
where $-1/(1+\beta^2)-1/(1+\beta)^2=\lambda-1$
 and 
\begin{equation}
\begin{bmatrix}
2&1\\
0&1 
\end{bmatrix}
\begin{bmatrix}
 \lambda-1&1\\
 1&0
\end{bmatrix}
\begin{bmatrix}
2&0\\
1&1 
\end{bmatrix}
=
\begin{bmatrix}
4\lambda&2\\
2&0 
\end{bmatrix}\eqqcolon Q,
\end{equation}
by noting that
\begin{equation}
Q_{11}\norm{x-y}^2
+2Q_{12}\innp{x-y,u-v}
+Q_{22}\norm{u-v}^2
\ge 0\quad  \forall (x,u), (y,v) \in \gra(R_{A}).
\end{equation} 
\qed

\section{Linear rates of convergence: three cases}
\label{sec:Motivation}

The Douglas--Rachford splitting algorithm, 
introduced by Lions and
Mercier~\cite{L-M79}, is a fundamental algorithm for solving 
monotone inclusion problems of the form
\begin{equation}
 \label{eq:P:sum}
 \text{Find $x\in X$ such that $0\in Ax+Bx$,}
\end{equation}
where  %$X$ is a real Hilbert space, and 
$A\colon X\rras X$ and $B\colon X\rras X$ are maximally 
monotone operators.   
%(See Section~\ref{sec:RA:prop} for a review of these
%and other definitions used in the paper.)
The algorithm is based on the iteration
\begin{equation}
 u_{n+1} = T_{\DR}u_n
 = \tfrac{1}{2} (\Id + R_BR_A)u_n,
\end{equation}
starting at arbitrary $u_0\in X$,  
where $R_A$ and $R_B$ are the 
{reflected resolvents} of $A$ and $B$.
%(defined as $R_A = 2J_A - I$ and $R_B = 2J_B-I$  where  the operators
%$J_A = (\Id + A)^{-1}$ and $J_B = (\Id + B)^{-1}$ are the 
%{resolvents}).
If the inclusion problem~\cref{eq:P:sum} has a
solution, then the iterates of $(u_k)_{k\in \NN}$ 
can be shown to converge weakly to
some point $u \in X$, where $u=T_{\DR}u$
and $x=J_Au$ solves~\cref{eq:P:sum}, see, e.g., 
\cite{BC2017,Comb04,Svaiter}.

In this section we review the results from \cite{Gis17}
on contraction properties of the Douglas--Rachford operator 
$T_{\DR} = (1/2)(\Id + R_BR_A)$.  These results are summarized in 
\cref{cor:T:prox:it}.
The following lemma shows that the three cases in
\cref{cor:T:prox:it} all have in common that $T_{\DR}$ is
the resolvent of a strongly monotone operator 
(hence, in view of 
\cref{prop:corres:A:RA}\ref{prop:corres:A:RA:vii}, 
a contraction).
We will see that this is a key difference with the new result
in Section~\ref{sec:Main}.

\begin{lemma}
\label{lem:T:prox:it}
Let $T_1\colon X\to X$ and $T_2\colon X\to X$ be nonexpansive. 
Define $T= \tfrac{1}{2}(\Id+T_2T_1)$ and $C= T^{-1}-\Id$.
Let $\alpha \in \left] 0,1\right[$.
Consider the following statements.
\begin{enumerate}[{\rm (a)}]
\item
\label{lem:T:prox:it:iv:a}
$C$ is $((1-\alpha)/\alpha)$-strongly monotone.
\item
\label{lem:T:prox:it:iv:c}
$-T_2T_1$ is $\alpha$-averaged\footnote{This is 
also known as negative averagedness of the operator $T_2T_1$
 \cite{Gis17}.}. 
%$\siff$
\item
\label{lem:T:prox:it:iv:d}
$T$ is $(1/\alpha)$-cocoercive.
%$\RA$ 
\item
\label{lem:T:prox:it:iv:e}
$T$ 
is a Banach contraction
with a constant
$\alpha$.
\end{enumerate}
Then 
\ref{lem:T:prox:it:iv:a} 
$\siff$
\ref{lem:T:prox:it:iv:c}
 $\siff$
\ref{lem:T:prox:it:iv:d} 
$\RA$
\ref{lem:T:prox:it:iv:e}. 
%\end{enumerate}
\end{lemma}
{\it Proof}
We first note that $T=J_C$ and $T_2T_1 = R_C$,
by definition of $T$ and $C$.
Hence $T$ is firmly nonexpansive 
by \cite[Theorem~2.1]{GK90}
and $C$ is maximally monotone by \cite[Theorem~2]{EckBer}.

\ref{lem:T:prox:it:iv:a} 
$\siff$
\ref{lem:T:prox:it:iv:c}
$\siff$
\ref{lem:T:prox:it:iv:d}: This follows from 
\cref{prop:corres:A:RA}\ref{prop:corres:A:RA:vii}
applied with 
$A$ replaced by $C$
and
$\mu$ replaced by $(1-\alpha)/\alpha$
and the fact that $T_2T_1 = R_C$
(see also \cite[Proposition~5.4]{Gis17}).
\ref{lem:T:prox:it:iv:d} 
$\RA$ \ref{lem:T:prox:it:iv:e}: 
This follows from
\cite[Proposition~23.13]{BC2017} (see also the comment after 
\cref{eq:Lip:coco}).
\qed

Reference \cite[Sections~5,~6~\&~7]{Gis17} contains a
comprehensive analysis of the rates of linear convergence of
the Douglas--Rachford method with optimal relaxation parameters
and step lengths, for the three cases presented in the next corollary. 
The key idea is that in each case, the Douglas--Rachford operator
is a contraction, as summarized below.

\begin{corollary}
\label{cor:T:prox:it} 
Let $\beta\geq \mu>0$.
Suppose that one of the following properties is satisfied.
\begin{enumerate}[{\rm(a)}]
\item
\label{cor:T:prox:it:A} 
$A$ is $(1/\beta)$-cocoercive and 
$B$ is $\mu$-strongly monotone. 
\item
\label{cor:T:prox:it:B}
$A$ is $(1/\beta)$-cocoercive and 
$\mu$-strongly monotone. 
\item 
\label{cor:T:prox:it:C}
$A$ is $\beta$-Lipschitz continuous and
$\mu$-strongly monotone.
\end{enumerate}
Then 
\begin{enumerate}
\item
\label{cor:T:prox:it:i} 
$-R_BR_A$ is $\alpha$-averaged for some 
$\alpha\in \left]0,1\right[$.
\item
\label{cor:T:prox:it:ii} 
$T_{\DR} = (1/2)(\Id + R_BR_A)$ 
is a Banach contraction with a contraction factor $\kappa\in \left]0,1\right[$. 
\end{enumerate}
The expressions for $\alpha$ and $\kappa$ are as follows.
\begin{equation}
\begin{array}{ll}
\mbox{Case (a):} 
& \displaystyle
\alpha=\kappa=\frac{1+\mu\beta}{1+\mu+\mu\beta} \\*[1ex]
\mbox{Case (b):} &\displaystyle
\alpha=\kappa
=\frac{1}{2}+\frac{1}{2}\left(\frac{1-2\mu+\mu\beta}{1+2\mu+\mu\beta}
\right)^{1/2} \\*[1ex]
\mbox{Case (c):} &\displaystyle
\alpha=\kappa
=\frac{1}{2}+\frac{1}{2}\left(\frac{1-2\mu+\beta^2}{1+2\mu+\beta^2}
\right)^{1/2}.
\end{array}
\end{equation}

\end{corollary}
{\it Proof}
We first discuss \ref{cor:T:prox:it:i}.
From \cite[Proposition~5.5]{Gis17},
Assumption \ref{cor:T:prox:it:A} implies that
$-R_BR_A$ is $\alpha$-averaged.
If Assumption \ref{cor:T:prox:it:B} holds,
then $R_A$ is a Banach contraction with factor 
\[
\kappa_1=\left(\frac{1-2\mu+\mu\beta}{1+2\mu+\mu\beta}
\right)^{1/2}
\]
(see \cite[Theorem~6.3]{Gis17} 
or \cref{prop:corres:A:RA}\ref{prop:corres:A:RA:ii}).
If Assumption \ref{cor:T:prox:it:C} holds, then $R_A$
is a Banach contraction with factor
\[
\kappa_2=\left(\frac{1-2\mu+\beta^2}{1+2\mu+\beta^2}
\right)^{1/2}
\]
(see \cite[Theorem~7.2]{Gis17} or 
\cref{prop:corres:A:RA}\ref{prop:corres:A:RA:iii}).
In both cases (\ref{cor:T:prox:it:B} and~\ref{cor:T:prox:it:C}), 
this implies 
that the compositions $R_BR_A$ 
and $-R_BR_A$ are Banach contractions 
with factors $\kappa_1$,
and  $\kappa_2$ respectively.  Hence 
$-R_BR_A$ is $((\kappa_1+1)/2)$-averaged 
(respectively $((\kappa_2+1)/2)$-averaged)
by \cite[Proposition~4.38]{BC2017}.

The second part \ref{cor:T:prox:it:ii} is proved by
combining \ref{cor:T:prox:it:i} and 
\cref{lem:T:prox:it}
applied with $T_1 = R_A$ and $T_2=R_B$,
and using the triangle inequality.
\qed

\section{Main results}
\label{sec:Main}
We now consider the Douglas--Rachford iteration under the following
assumptions:
 
 \begin{equation}
\label{eq:nonexp:assmp}
A\colon  X\to X \text{~is $\beta$-Lipschitz continuous
 and monotone, and $\beta>0$}
\end{equation}
 and that 
  \begin{equation}
\label{eq:stmnono:assmp}
 B\colon X\rras X
 \text{~is maximally monotone and 
 $\mu$-strongly monotone, and $\mu>0$}. 
\end{equation}
This case is not covered by \cref{cor:T:prox:it},
 and is significantly
different in nature, because these two properties 
in \cref{eq:nonexp:assmp} 
and \cref{eq:stmnono:assmp}
do not imply that
$-R_BR_A$ is averaged, 
as shown by the following example.

\begin{example}
\label{ex:nonexp:notav}
Suppose that $X=\RR^2$ and define
\begin{equation}
A= \begin{bmatrix}
0&1\\
-1&0
\end{bmatrix}, \quad
B= N_{\{0\}}.
\end{equation}
Then
$A$ is monotone and nonexpansive (hence $1$-Lipschitz continuous),
$B$ is maximally monotone 
and $\mu$-strongly monotone for every $\mu>0$,
$R_A=-A$ and $R_B=-\Id$.
Hence, $-R_BR_A=-A$ which is \emph{not} averaged.
\end{example}

The main results in this section are
\cref{thm:DR:stmono}
and \cref{thm:DR:nonex:stmono:linear} below.
We first prove a more general result on an averaged composition
of a $\beta$-Lipschitz continuous operator 
 and an averaged operator.
 
\begin{proposition}
\label{prop:nonexp:gen:or:DR}
Let $R\colon X\to X$ be such that $-R$ is
 $\alpha$-averaged, with $\alpha\in \left[0,1\right[$.
Let $M\colon X \to X$ be nonexpansive such that
$(\forall (x,y)\in X\times X)$ 
\begin{equation}
\label{eq:assump:Lips:RA}
\innp{x-y,Mx-My}\ge -\lambda \norm{x-y}^2,\text{~with~}
\lambda\in \left[-1, 1\right[.
\end{equation}
Define 
\begin{equation}
T=\tfrac{1}{2}(\Id+RM), \qquad
\widetilde{T}=\tfrac{1}{2}(\Id+MR).
\end{equation}
Then the following hold:
\begin{enumerate}
\item
\label{prop:nonexp:gen:or:DR:i}
$\Id+(\alpha-1)M$ is Lipschitz continuous with constant
$\sqrt{1+(1-\alpha)^2+2\lambda (1-\alpha)}<2-\alpha<2$.
\item
\label{prop:nonexp:gen:or:DR:ii}
$T$ is Lipschitz continuous 
 with constant
\begin{equation}
\label{eq:rate:abstract}
%\tfrac{1}{2}(\sqrt{(1-\alpha)^2+1}+\alpha)
 \tfrac{1}{2}\Big(\sqrt{1+(1-\alpha)^2+2\lambda (1-\alpha)}+\alpha\Big)<1.
 \end{equation} 
Hence, $T$ is a Banach contraction
 and $\fix T$ is a singleton.
\item
\label{prop:nonexp:gen:or:DR:iii}
If $M$ is linear, then
$\widetilde{T}$ is Lipschitz continuous 
 with constant given in \cref{eq:rate:abstract}.
Hence, $\widetilde{T}$ is a Banach contraction
 and $\fix \widetilde{T}$ is a singleton.
\end{enumerate}
\end{proposition}
{\it Proof}
\ref{prop:nonexp:gen:or:DR:i}:
Set $S= \Id+(\alpha-1)M$
 and let $(x,y)\in X\times X$.
Then
 \begin{subequations} 
  \begin{align}
\norm{Sx-Sy}^2
&=\norm{x-y}^2+(1-\alpha)^2\norm{Mx-My}^2
-2(1-\alpha)\innp{x-y,Mx-My}\\
&\le \norm{x-y}^2+(1-\alpha)^2\norm{Mx-My}^2
+2\lambda (1-\alpha)\norm{x-y}^2\\
& \le (1+(1-\alpha)^2+2\lambda (1-\alpha))\norm{x-y}^2.
\end{align}
\end{subequations}
The first inequality follows from \cref{eq:assump:Lips:RA}
 and the 
second inequality follows from the nonexpansiveness of $M$. 
Finally note that, because 
$-1\le \lambda<1$, we learn that
$\sqrt{1+(1-\alpha)^2+2\lambda(1-\alpha)}
<\sqrt{1+(1-\alpha)^2+2(1-\alpha)}
=\sqrt{1+1-2\alpha+\alpha^2+2-2\alpha}
=\sqrt{(2-\alpha)^2}
=2-\alpha<2$.

\ref{prop:nonexp:gen:or:DR:ii}:
Since $-R$ is $\alpha$-averaged,
we have $R=(\alpha-1)\Id+\alpha N$ for some nonexpansive
$N\colon X\to X$.
Substituting this in the definition of $T$, we get
$T=\tfrac{1}{2}(\Id+(\alpha-1)M+\alpha NM)$.
It follows from the triangle inequality, 
\ref{prop:nonexp:gen:or:DR:i},
and the nonexpansiveness of $M$ and $N$ that $T$ is Lipschitz 
continuous with a constant 
\begin{equation}
 \tfrac{1}{2}\Big(\sqrt{1+(1-\alpha)^2+2\lambda 
  (1-\alpha)}+\alpha\Big)
<\tfrac{1}{2}(2-\alpha+\alpha)=1.
\end{equation}
   
\ref{prop:nonexp:gen:or:DR:iii}:
As in \ref{prop:nonexp:gen:or:DR:ii},
we write $R$ as $R=(\alpha-1)\Id+\alpha N$ with $N$ nonexpansive.
Then 
\[
\widetilde{T} 
=\tfrac{1}{2} (\Id + M((\alpha-1)\Id + \alpha N))
=\tfrac{1}{2}(\Id+(\alpha-1)M+\alpha MN).
\]
The second identity follows from linearity of $M$.
Now the proof of \ref{prop:nonexp:gen:or:DR:iii}
is similar to \ref{prop:nonexp:gen:or:DR:ii}.
\qed

We are now ready for our main results.

\begin{theorem}
\label{thm:DR:stmono}
Suppose that 
$A\colon X\to X$ is monotone and %linear, skew and 
$\beta$-Lipschitz continuous with $\beta>0$,
%nonexpansive,
and that
$B\colon X \rras X$ is maximally monotone and 
$\mu $-strongly monotone with $\mu > 0$.
Let $x_0\in X$,
%let \hl{$\gamma>0$},
let $T=\tfrac{1}{2}\Big(\Id+R_{ B}R_{ A}\Big)$.
Then the following hold:
 \begin{enumerate}
 \item
 \label{thm:DR:stmono:i}
 $(x_n)_\nnn = (T^nx_0)_\nnn$ converges strongly 
 to some $\overline{x} \in X$, with a linear rate $r$, where
\begin{equation}
\label{eq:rate:concrete}
r=
 \tfrac{1}{2(1+\mu)}
 \Big(\sqrt{2\mu^2+2\mu+1
 +2\left(1-\tfrac{1}{(1+\beta)^2}
 -\tfrac{1}{1+\beta^2}\right)\mu(1+\mu)}
 +1\Big)<1.
\end{equation} 
 \item
  \label{thm:DR:stmono:ii}
  $(J_{ A}x_n)_\nnn$ converges strongly to 
  $J_{ A}\overline{x}$ with a linear rate 
  $r$ given in \cref{eq:rate:concrete}.
   \end{enumerate}
Moreover, $\fix R_{B}R_{ A}
=\fix T=\{\overline{x}\}$,
%$\fix R_{\gamma A}R_{\gamma B}
%=\fix \widetilde{T}=\{\hat{x}\}$
and 
$\zer(A+B)=\{J_{ A}\overline{x}\}$.
%=\{J_{\gamma B}\hat{x}\}$.
\end{theorem}
{\it Proof}
Since $ A$ is 
monotone and 
$ \beta $-Lipschitz continuous,
we have
$R_{ A}$ is nonexpansive
 and $(\forall (x,y)\in X\times X)$
  \begin{equation}
 \innp{x-y, R_{ A} x-R_{ A}y}\ge -
\left(1-\tfrac{1}{(1+\beta)^2}
 -\tfrac{1}{1+\beta^2}\right)\normsq{x-y},
 \end{equation}
by \cref{lem:RA:Lips:beta}\ref{lem:RA:Lips:beta:ii}.
Since $ B$ is $\mu $-strongly monotone, 
$-R_{ B}$ is $( 1+\mu)^{-1}$-averaged 
(see 
\cite[Proposition~5.4]{Gis17} or 
\cref{prop:corres:A:RA}\ref{prop:corres:A:RA:vii}).
\ref{thm:DR:skew:stmono:i}:
The claim of strong convergence follows from
\cite[Theorem~26.11(vi)(a)]{BC2017}.
The rate $r$ follows from 
\cref{prop:nonexp:gen:or:DR}\ref{prop:nonexp:gen:or:DR:ii}
applied with $\alpha = (1+\mu)^{-1}$, 
$\lambda=(1-({1}/{(1+\beta)^2})
 -({1}/{1+\beta^2}))$,
$M=R_{ A}$, and $R=R_{ B}$.
\ref{thm:DR:skew:stmono:ii}:
This is a direct consequence of 
\ref{thm:DR:skew:stmono:i}
 and the fact that $J_{ A}$ is
 (firmly) nonexpansive.
\qed

When $A$ is linear, similar conclusion to
that of  \cref{thm:DR:stmono}
holds if we switch
the order of the operators in the Douglas--Rachford iteration.
\begin{theorem}
\label{thm:DR:nonex:stmono:linear}
Suppose that 
$A\colon X\to X$ is monotone, $\beta$-Lipschitz 
continuous with $\beta >0$ and linear,
and that
$B\colon X \rras X$ is maximally monotone and 
$\mu $-strongly monotone with $\mu > 0$.
Let $\widetilde{x}_0\in X$,
%let \hl{$\gamma>0$} 
and 
let 
$\widetilde{T}=\tfrac{1}{2}\Big(\Id+R_{ A}R_{ B}\Big)$.
Then the following hold:
\begin{enumerate}
 \item
 \label{thm:DR:nonex:stmono:iii}
 $(\widetilde{x}_n)_\nnn = (\widetilde{T}^n x_0)_\nnn$
 converges strongly to some $\hat{x}$ 
 with a linear rate   $r$ given in \cref{eq:rate:concrete}.
 \item
 \label{thm:DR:nonex:stmono:iv}
  $(J_{ B}\widetilde{x_n})_\nnn$ converges strongly 
  to $J_{ B}\hat{x}$ with a linear rate 
  $r$ given in \cref{eq:rate:concrete}.

 \end{enumerate}
Moreover, $\fix R_{ A}R_{ B}
=\fix \widetilde{T}=\{\hat{x}\}$
  and $\zer (A+B)=\big\{J_{ B}\hat{x}\big\}$.
\end{theorem}
{\it Proof}
Proceed as in the proof of
\cref{thm:DR:stmono}\ref{thm:DR:stmono:i}--\ref{thm:DR:stmono:ii},
but use 
\cref{prop:nonexp:gen:or:DR}\ref{prop:nonexp:gen:or:DR:iii}
and the fact that $R_{ A}$ is linear.
\qed

We conclude this section with the following remark.
\begin{remark}\
It is not clear whether or not the conclusion of 
\cref{prop:nonexp:gen:or:DR}\ref{prop:nonexp:gen:or:DR:iii}
remains true if we drop the assumption of linearity.
Any counterexample to show failure of the conclusion 
in the absence of linearity must feature 
nonexpansive 
operator that
satisfies \cref{eq:assump:Lips:RA}
which is \emph{neither} linear \emph{nor} averaged,
because
if $M$ is averaged we have that $\widetilde{T}$ is a Banach contraction
by \cite[Proposition~3.9]{Gis17}. 
Note, however, that the result for linear $M$ covers
important applications, such as the
the primal-dual Douglas--Rachford method 
discussed in Section~\ref{sec:Applications}.

\end{remark}

\section{The linear skew case and
application to primal-dual Douglas--Rachford method}
\label{sec:Applications}
The main goal of this section is to
 prove linear convergence of 
the primal-dual 
Douglas--Rachford method
discussed in
 \cite[Sections~3.1\&3.2]{OV14}
 (see also  
 \cite{B-AC11} for a more general framework)
 when applied to solve 
the monotone inclusion
\cref{eq:mono:inc:bAbB} below,
under additional assumptions on the
underlying operators.

In the following we show that when $A$
is linear and skew then the rate in 
\cref{eq:rate:concrete}
%\cref{sec:Main}
is improved.
We first start with the following lemma
which shows that when
$A$ is linear and skew,
the bounds in \cref{lem:RA:Lips:beta}
can be tightened.
\begin{lemma}
\label{lem:RA:skewop}
Suppose that $A\colon X\to X$ is linear, skew, i.e., $A=-A^*$ and
$\beta$-Lipschitz continuous with $\beta>0$. 
Let $x\in X$.
Then the following hold:
%\begin{equation}\end{equation}
\begin{enumerate}
\item
\label{lem:RA:skewop:i}
$R_A$ is an isometry, i.e., $\norm{R_Ax}=\norm{x}$.
%Hence $R_A^{-1}=R_A^*$.
\item
\label{lem:RA:skewop:ii}
$\norm{x}^2\le (1+\beta^2)\norm{J_A x}^2$.
\item
\label{lem:RA:skewop:iii}
$J_A$ is $\tfrac{1}{\beta^2+1}$-strongly monotone.
\item
\label{lem:RA:skewop:iv}
$\innp{x,R_Ax}\ge \big(\tfrac{2}{1+\beta^2}-1\big)\norm{x}^2$.
\end{enumerate}
\end{lemma}
{\it Proof}
Set $u=J_Ax$ and note that 
$Au=x-u$.
Now, since $A$ is skew, in view of
 \cref{eq:A:JA} we have
 \begin{equation}
 \label{eq:skew:gr}
% (\forall x\in X)\quad 
 \innp{u,x-u}
 =\innp{u,Au}=0.
 \end{equation}
\ref{lem:RA:skewop:i}:
Using \cref{eq:skew:gr} we have 
\begin{subequations}
\begin{align}
\norm{R_A x}^2
&=\normsq{J_A x-J_{A^{-1}}x}=\normsq{u-Au}\\
&=\normsq{u}-2\innp{u,Au}+\normsq{Au}\\
&=\normsq{u}+2\innp{u,Au}+\normsq{Au}\\
&=\normsq{u+Au}=\normsq{x}.
\end{align}
\end{subequations}
\ref{lem:RA:skewop:ii}:
Indeed, using \cref{eq:skew:gr} we have 
\begin{subequations}
\begin{align}
\normsq{x}
&=\normsq{u}+2\innp{u,Au}+\normsq{Au}
=\normsq{u}+\normsq{Au}\\
&\le \normsq{u}+\beta^2\normsq{u}
=(1+\beta^2) \normsq{u},
\end{align}
\end{subequations}
where the inequality follows from the 
$\beta$-Lipschitz continuity of $A$.
% and the fact that $(J_A x,J_{A^{-1}}x)
%\in \gra A$.

\ref{lem:RA:skewop:iii}:
It follows from 
\cref{lem:RA:Lips:beta}\ref{lem:RA:Lips:beta:iv}
that
$(1+\beta^2)\normsq{x-J_Ax}\le \beta^2\normsq{x}$.
Expanding yields 
$(1+\beta^2)(\normsq{x}+\normsq{J_Ax}-2\innp{x,J_Ax})\le \beta^2\normsq{x}$.
Equivalently,
$2(1+\beta^2)\innp{x,J_Ax}\ge \normsq{x}+(1+\beta^2)\normsq{J_{A}x}$.
Now combine with \ref{lem:RA:skewop:ii}.

\ref{lem:RA:skewop:iv}:
We have 
$\innp{x,R_A x}
=\innp{x, 2J_A x-x}
=2\innp{x,J_Ax }-\normsq{x}
\ge \big(\tfrac{2}{\beta^2+1}-1\big)\normsq{x}$,
where the inequality follows from
\ref{lem:RA:skewop:iii}.
\qed

\begin{theorem}
\label{thm:DR:skew:stmono}
Suppose that 
$A\colon X\to X$ is linear, skew (hence monotone) and 
$\beta$-Lipschitz continuous with $\beta>0$,
%nonexpansive,
and that
$B\colon X \rras X$ is maximally monotone and 
$\mu $-strongly monotone with $\mu > 0$.
Let $x_0\in X$,
%let $\widetilde{x}_0\in X$,
%let \hl{$\gamma>0$},
let $T=\tfrac{1}{2}\big(\Id+R_{ B}R_{ A}\big)$
and let 
$\widetilde{T}=\tfrac{1}{2}\big(\Id+R_{ A}R_{ B}\big)$.
Then the following hold:
 \begin{enumerate}
 \item
 \label{thm:DR:skew:stmono:i}
 $(x_n)_\nnn = (T^nx_0)_\nnn$ converges strongly 
 to some $\overline{x} \in X$, with a linear rate $r$, where
\begin{equation}
\label{eq:rate:concrete:skew}
r(\beta,\mu)=
 \tfrac{1}{2(1+\mu)}
 \Big(\sqrt{2\mu^2+2\mu+1
 +2\left(1-\tfrac{2}{1+\beta^2}\right)\mu(1+\mu)}
 +1\Big)<1.
\end{equation} 
 \item
  \label{thm:DR:skew:stmono:ii}
  $(J_{ A}x_n)_\nnn$ converges strongly to 
  $J_{ A}\overline{x}$ with a linear rate 
  $r$ given in \cref{eq:rate:concrete:skew}.
   \item
 \label{thm:DR:skew:stmono:iii}
 $(\widetilde{x}_n)_\nnn = (\widetilde{T}^n x_0)_\nnn$
 converges strongly to some $\hat{x}$ 
 with a linear rate   $r$ given in \cref{eq:rate:concrete:skew}.
 \item
 \label{thm:DR:skew:stmono:iv}
  $(J_{ B}\widetilde{x_n})_\nnn$ converges strongly 
  to $J_{ B}\hat{x}$ with a linear rate 
  $r$ given in \cref{eq:rate:concrete:skew}.

   \end{enumerate}
Moreover, $\fix R_{ B}R_{ A}
=\fix T=\{\overline{x}\}$,
$\fix R_{ A}R_{ B}
=\fix \widetilde{T}=\{\hat{x}\}$
and 
$\zer(A+B)=\{J_{ A}\overline{x}\}
=\{J_{ B}\hat{x}\}$.
\end{theorem}
{\it Proof}
Proceed as in the proof of 
\cref{thm:DR:stmono}
for 
\ref{thm:DR:skew:stmono:i}--\ref{thm:DR:skew:stmono:ii}
(respectively
\cref{thm:DR:nonex:stmono:linear}
for 
\ref{thm:DR:skew:stmono:iii}--\ref{thm:DR:skew:stmono:iv})
in view of \cref{lem:RA:skewop}\ref{lem:RA:skewop:iv}.
\qed

The contraction factor in \cref{eq:rate:concrete:skew}
is sharp as we illustrate in \cref{ex:sharp:rate} below.
\begin{example}[\thmtit{sharpness of the
contraction factor}]
\label{ex:sharp:rate}
Let $\beta>0$ and let $\mu>0$.
Suppose that
$X=\RR^2$,
\begin{equation}
A=\beta\begin{bmatrix}
0&1\\
-1&0
\end{bmatrix}, \qquad
B=\mu\Id+N_{\{0\}\times \RR}.
\end{equation}
Then 
$A$ is $\beta$-Lipschitz continuous
and monotone, 
$B$ is $\mu$-strongly monotone and
%for every $\gamma>0$
%the reflected resolvents
%of $ A$ and $ B$
%are
\begin{equation} \label{e-RB}
R_{ A}=
\begin{bmatrix}
\tfrac{2}{\beta^2+1}-1
&-\tfrac{2\beta}{\beta^2+1}\\
\tfrac{2\beta}{\beta^2+1}
&\tfrac{2}{\beta^2+1}-1
\end{bmatrix},
\qquad
R_{ B}=
\begin{bmatrix}
-1&0\\
0&\frac{1- \mu}{1+\mu}
\end{bmatrix}.
\end{equation}
Therefore,
\begin{equation}
T=\tfrac{1}{2}\left(\Id+R_{ B}R_{ A}\right)
=\frac{1}{\beta^2+1}
\begin{bmatrix}
\beta^2
& \beta\\
\frac{ \beta(1-\mu)}{1+\mu}
&\frac{1+\beta^2\mu}{1+\mu}
\end{bmatrix}, 
\end{equation}
\begin{equation}
\norm{T}=
 \tfrac{1}{2(1+\mu)}
 \Big(\sqrt{2\mu^2+2\mu+1
 +2\left(1-\tfrac{2}{1+\beta^2}\right)
 \mu(1+\mu)}
 +1\Big).
%\frac{1}{2(1+\mu)}
%\left(\sqrt{2\mu^2+2\mu+1}+1\right).
\end{equation}
\end{example}
{\it Proof}
The claim about $R_{ A}$
%follows from \cref{eq:cal:RA:skew}.
is straightforward to verify.
By \cite[Example~23.4~and~Corollary~3.24(iii)]{BC2017}, 
we have
\begin{subequations}
\begin{align}
J_{B}& %=(\Id+B)^{-1}
=((1+\mu)\Id+N_{\{0\}\times \RR})^{-1} \\
& =((1+\mu)(\Id+N_{\{0\}\times \RR}))^{-1}\\
&=J_{N_{\{0\}\times \RR}}\circ \tfrac{1}{1+\mu}\Id \\
& =\tfrac{1}{1+\mu}P_{\{0\}\times \RR} \\
&=\tfrac{1}{1+\mu}
\begin{bmatrix}
0&0\\
0&1
\end{bmatrix}
.
\end{align}
\end{subequations}
The expression for $R_{ B}$ in~\cref{e-RB} 
and the formula for $T$ readily follows.
A routine calculation yields that the eigenvalues of  $TT\tran$
are
\begin{equation}
%\resizebox{0.9\hsize}{!}{ $
\frac{2\beta^2\mu^2
+2\beta^2\mu+\beta^2+1
\pm\sqrt{
(2\beta^2\mu^2
+2\beta^2\mu+\beta^2+1)^2
-4(1+\mu)^2(1+\beta^2)\beta^4\mu^2
}}{2(1+\mu)^2(1+\beta^2)}.
%$}
\end{equation}
Hence, 
\begin{subequations}
\begin{align}
\norm{T} &=
\norm{TT\tran}^{1/2}
%=\tfrac{1}{2(1+\mu)^2}
%\left(\mu^2+\mu+1+\sqrt{2\mu^2+2\mu+1}\right)
\\
&=
\resizebox{0.8\hsize}{!}{ $
\sqrt{\frac{2\beta^2\mu^2
+2\beta^2\mu+\beta^2+1
+\sqrt{
(2\beta^2\mu^2
+2\beta^2\mu+\beta^2+1)^2
-4(1+\mu)^2(1+\beta^2)\beta^4\mu^2
}}{2(1+\mu)^2(1+\beta^2)}}
$}
\\
&= \frac{1}{2(1+\mu)}
 \Big(\sqrt{2\mu^2+2\mu+1
 +2\left(1-\tfrac{2}{1+\beta^2}\right)
 \mu(1+\mu)}
 +1\Big).
\end{align}
\end{subequations}
\qed

\begin{figure}
%\begin{minipage}{1.2\textwidth}
\centering
\begin{overpic}[scale=0.24]{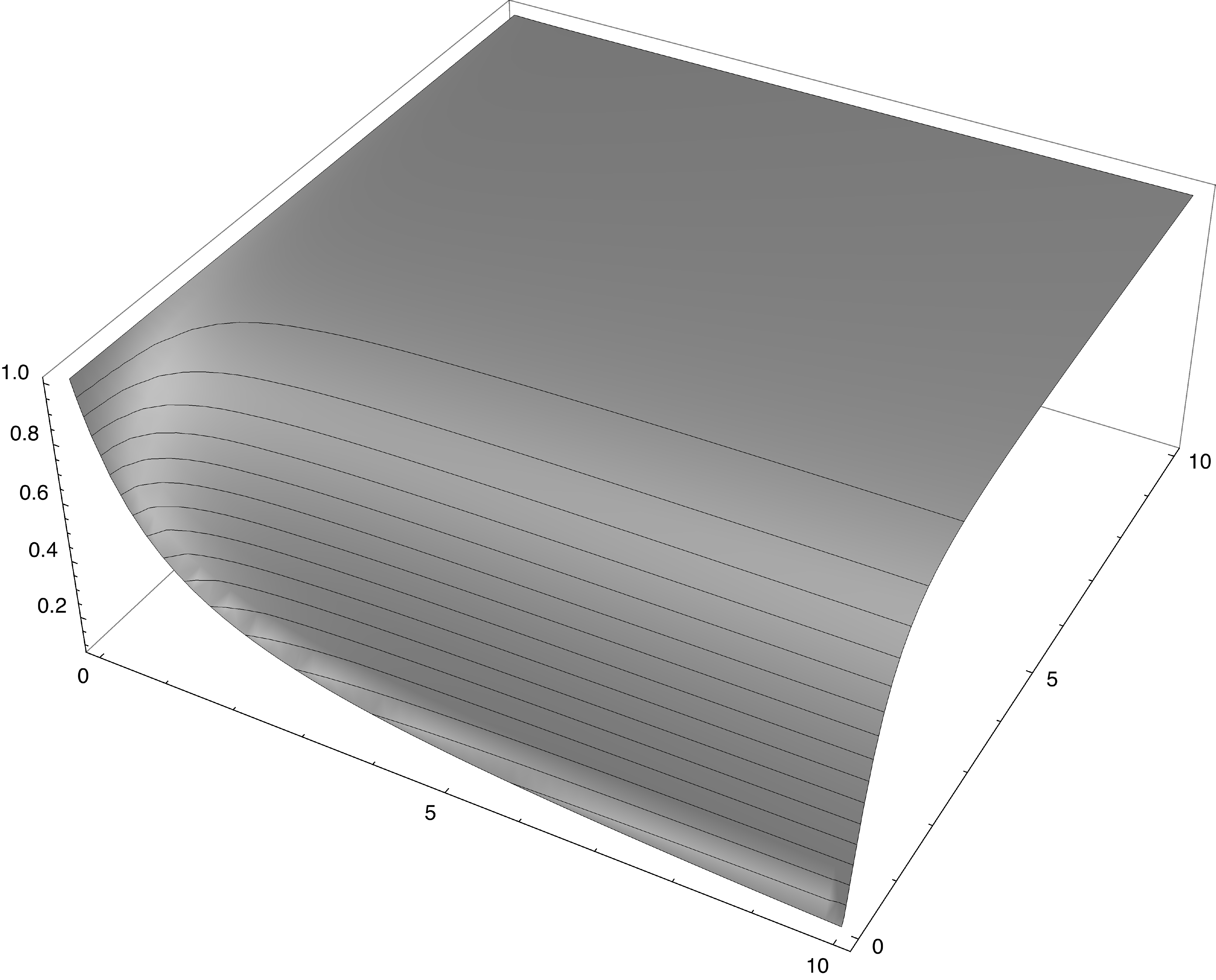}
%\begin{overpic}[width=0.5\textwidth,grid,tics=10]{plot}
 \put (30,10) {\scriptsize$\beta$}
  \put (90,22) {\scriptsize$\mu$}
  \put (-5,28) {\rotatebox{102}{\scriptsize${r(\beta,\mu)}$}}
\end{overpic}
%\includegraphics[scale=0.6]{cbar}
%}
%%    \caption{Pijl}
%\end{minipage}
\caption{A \texttt{Mathematica} \cite{Wolfram} snapshot.
Shown is the rate $r=r(\beta,\mu)$
given in
\cref{eq:rate:concrete:skew} for the case
when $A$ is $\beta$-Lipschitz continuous
and monotone and $\beta>0$ and 
$B$ is $\mu$-strongly monotone. 
}
\label{fig3}
\end{figure}

\begin{figure}
\begin{center}
%\begin{minipage}{.7\textwidth}
\begin{tikzpicture}
\begin{axis}[
       width=0.45\textwidth,
       height=0.4\textwidth,
    axis lines = left,
     xmin=0, xmax=5, xstep=2,
    ymin=0, ymax=1,
    legend pos=south east,
    xlabel = $\mu$,
    ylabel = {rate of convergence},
         legend columns=2, 
     legend style={at={(1,0.01)},
     font=\fontsize{8}{9}\selectfont,
     /tikz/column 2/.style={
                column sep=4pt,
            },},
]
\addplot [
domain=0:10,
    samples=300,
    color=black,
    loosely dotted,
very thick
    ]
    {(1+sqrt(1+2*x+2*x^2+2*x*(1+x)*(1-2*(1/(1+(0.2)^2)))))/(2+2*x)};
\addlegendentry{$\beta=0.2$}
\addplot [
domain=0:10,
    samples=300,
    color=minhblue,
 thick
    ]
    {(1+sqrt(1+2*x+2*x^2+2*x*(1+x)*(1-2*(1/(1+(0.5)^2)))))/(2+2*x)};
\addlegendentry{$\beta=0.5$}
\addplot [
domain=0:10,
    samples=300,
    color=Brown, %Turquoise, %RoyalPurple,
    densely dotted,
 thick
    ]
    {(1+sqrt(1+2*x+2*x^2+2*x*(1+x)*(1-2*(1/(1+1^2)))))/(2+2*x)};
\addlegendentry{$\beta=1$}
\addplot [
domain=0:10,
    samples=300,
    color=red,
 dashed,
    dash pattern=on 5.6pt off 3pt,
 thick
    ]
     {(1+sqrt(1+2*x+2*x^2+2*x*(1+x)*(1-2*(1/(1+2^2)))))/(2+2*x)};
    \addlegendentry{$\beta=2$}
    \addplot [
domain=0:10,
    samples=300,
    color=OliveGreen,
    dashdotted,
       dash pattern=on 6pt off 2pt on \the\pgflinewidth off 2pt,
 thick
    ]
     {(1+sqrt(1+2*x+2*x^2+2*x*(1+x)*(1-2*(1/(1+5^2)))))/(2+2*x)};
    \addlegendentry{$\beta=5$ }
%\addlegendentry{\cref{thm:DR:nonex:stmono}}
\end{axis}
\end{tikzpicture}
%\end{minipage}
%\begin{minipage}{.7\textwidth}
\hspace{1cm}
 \begin{tikzpicture}
\begin{axis}[
       width=0.45\textwidth,
       height=0.4\textwidth,
    axis lines = left,
     xmin=0, xmax=5, xstep=2,
    ymin=0, ymax=1,
    legend pos=south east,
    xlabel = $\beta$,
    ylabel = {rate of convergence},
    legend columns=2, 
     legend style={at={(1,0.01)},
     font=\fontsize{8}{9}\selectfont,
     /tikz/column 2/.style={
                column sep=4pt,
            },},
]
\addplot [
domain=0:10,
    samples=300,
    color=black,
    loosely dotted,
 very thick
    ]
    {(1+sqrt(1+2*(0.2)+2*(0.2)^2+2*(0.2)*(1+(0.2))*(1-2*(1/(1+(x)^2)))))/(2+2*(0.2))};
\addlegendentry{$\mu=0.2$}

\addplot [
domain=0:10,
    samples=300,
    color=minhblue,
 thick
    ]
    {(1+sqrt(1+2*(0.5)+2*(0.5)^2+2*(0.5)*(1+(0.5))*(1-2*(1/(1+(x)^2)))))/(2+2*(0.5))};
\addlegendentry{$\mu=0.5$}
\addplot [
domain=0:10,
    samples=300,
    color= Brown,%Turquoise, %RoyalPurple,
    densely dotted,
 thick
    ]
       {(1+sqrt(1+2*(1)+2*(1)^2+2*(1)*(1+(1))*(1-2*(1/(1+(x)^2)))))/(2+2*(1))};
\addlegendentry{$\mu=1$}
%\addplot [
%domain=0:10,
%    samples=300,
%    color=red,
% thick
%    ]
%    {x/(1+x)};
\addplot [
domain=0:10,
    samples=300,
    color=red,
    %loosely 
    dashed,
    dash pattern=on 5.6pt off 3pt,
 thick
    ]
       {(1+sqrt(1+2*(2)+2*(2)^2+2*(2)*(1+(2))*(1-2*(1/(1+(x)^2)))))/(2+2*(2))};
    \addlegendentry{$\mu=2$}
    \addplot [
domain=0:10,
    samples=300,
    color=OliveGreen,
    dashdotted,
           dash pattern=on 6pt off 2pt on \the\pgflinewidth off 2pt,
 thick
    ]
     {(1+sqrt(1+2*(5)+2*(5)^2+2*(5)*(1+(5))*(1-2*(1/(1+(x)^2)))))/(2+2*(5))};
    \addlegendentry{$\mu=5$ }
%\addlegendentry{\cref{thm:DR:nonex:stmono}}
\end{axis}
\end{tikzpicture}
%\end{minipage}
\end{center}
\caption[foo bar]{Left:
Shown are the optimal rates of convergence 
given in \cref{eq:rate:concrete:skew}
as functions of $\mu$ for 
$\beta=0.2$ (black loosely-dotted line),
$\beta=0.5$ (blue solid line),
$\beta=1$ (brown densely-dotted line), 
$\beta=2$ (red dashed line)
and $\beta=5$ (green dash-dotted line).
Right:
Shown are the optimal rates of convergence 
given in \cref{eq:rate:concrete:skew}
as functions of $\beta$ for 
$\mu=0.2$ (black loosely-dotted line),
$\mu=0.5$ (blue solid line),
$\mu=1$ (brown densely-dotted line), 
$\mu=2$ (red dashed line)
and $\mu=5$ (green dash-dotted line).
}
\label{fig:mu:beta}
\end{figure}
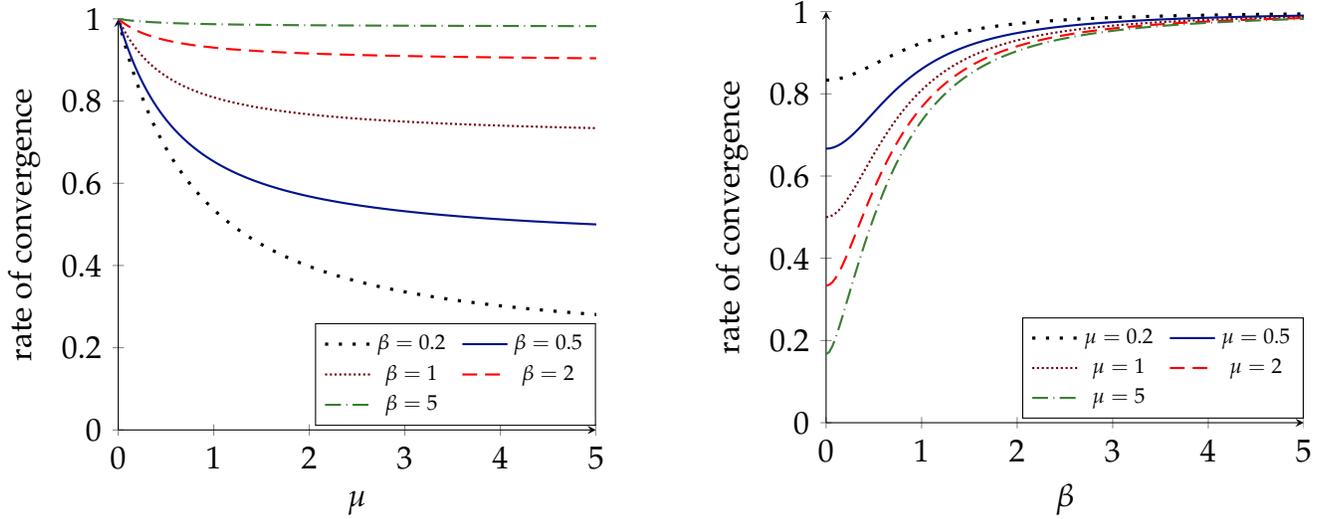

In Figure~\ref{fig3} we provide a
 plot of the rate in \cref{eq:rate:concrete:skew}
as a function of $\beta$ and $\mu$.
Figure~\ref{fig:mu:beta} 
provides plots of the rate in \cref{eq:rate:concrete:skew}
as a function of $\mu $ (respectively $\beta$)
for some concrete values of $\beta $ (respectively $\mu$).

\begin{remark}
Suppose that $\gamma>0$.
Then, in the setting of \cref{thm:DR:skew:stmono},
 the rate obtained when 
iterating $T=(1/2)(\Id+R_{\gamma B}R_{\gamma A})$ 
or 
$T=(1/2)(\Id+R_{\gamma A}R_{\gamma B})$
is $r(\gamma \beta, \gamma \mu)$,
where  $r$ is defined in \cref{eq:rate:concrete:skew}.  
However, 
unlike the rates presented in \cite{Gis17},
the rate given in \cref{eq:rate:concrete:skew}
cannot be easily optimized as a function of the
step-length $\gamma$.
Indeed, suppose that $\beta=\mu=1$. Then 
$r(\gamma \beta, \gamma \mu)
=h(\gamma)= \tfrac{1}{2(1+\gamma)}
 \Big(\sqrt{2\gamma^2+2\gamma+1
 +2\left(1-({2}/{(1+\gamma^2)})\right)\gamma(1+\gamma)}
 +1\Big)$. One can show that $h'(\gamma)=0$ reduces 
 to solving the quintic
 $4\gamma^5 +5\gamma^4 +12\gamma^3 +2\gamma^2 -3=0$.
 Of course we can solve numerically for the optimal value of 
 $\gamma$, which in this case yields $\gamma\approx0.4815$. 
 %$\gamma=.4815371581$
\end{remark}

Throughout the remainder of this section, we assume that
$Y$ is a real Hilbert space\footnote{A 
finite-dimensional example is
$(X,Y)=(\RR^n,\RR^m)$.},
that 
$L\colon X\to Y$ is nonzero and linear,
that\footnote{A \emph{closed} function is
 also known as 
lower semicontinuous.} 
%\begin{equation}%[box=\mybluebox]{equation}
%\label{f:assmp}
$f\colon X\to \left]-\infty,+\infty\right]$
is $\sigma$-strongly convex and closed,
%\end{equation} 
and that
$g\colon X\to \RR$
is convex and $\grad g$ is $\beta$-Lipschitz
continuous
for some $\beta>0$.

Consider the monotone inclusion:
\begin{equation}
\label{eq:mono:inc:bAbB}
\text{Find $(x,y)\in X\times Y$ such that
$0\in 
{\bf A}(x,y)+{\bf B}(x,y)$,}
\end{equation}
where\footnote{Here and 
elsewhere we  use $f^*$ to denote 
the \emph{convex conjugate} 
(this is also known as 
the \emph{Fenchel or Legendre conjugate})
of
$f$ defined at $u\in X$ 
as $f^*(u)
=\sup_{x\in X}\stb{\innp{u,x}-f(x)}$.} 
\begin{equation}
\label{eq:def:bA}
{\bf A}\colon X\times Y \to X\times Y\colon (x,y)
\mapsto (L^*y,-Lx),\quad 
%\text{~~and~~}
{\bf B}\colon X\times Y\rras X\times Y\colon(x,y)
 \mapsto \pt {f}(x)\times \pt {g}^*(y).
\end{equation}
One can check that
%\footnote{Let $(u,v)$ be a unit vector in 
%$X\times Y$ such that $0\not\in\{\norm{u}_X,\norm{v}_Y\}$
% and note that
%$\norm{\bA^*\bA(u,v)}^2
%=\norm{L^*Lu}_X^2+\norm{LL^*v}_Y^2
%=\norm{u}^2_X\norm{L^*L\tfrac{u}{\norm{u}}}_X^2
%+\norm{v}^2_Y\norm{LL^*\tfrac{v}{\norm{v}}}_Y^2$.
%Now use that $\norm{u}^2_X+\norm{v}^2_Y=1$,
%that $\norm{L}^2=\norm{LL^*}=\norm{L^*L}$
% and that $\norm{\bA}^2=\norm{\bA\bA^*}$ to learn that
%$\norm{\bA}\le \norm{L}$. The reverse inequality follows 
%similarly by noting that
%for a unit vector $u\in X$ we have
%$\norm{\bA\bA^*}^2\ge 
%\sup_{\norm{u}_X=1}\norm{\bA^*\bA(u,0)}^2
%=\sup_{\norm{u}_X=1} \norm{L^*Lu}_X^2
%=\norm{L^*Lu}_X^2$.
%}
 $\norm{{\bf A}}=\norm{L}\neq 0$. Hence,
\begin{equation}
\text{\label{lem:aux:pddr:b:i}
${\bf A}$ is Lipschitz continuous 
with the sharp constant $\norm{L}$.}
\end{equation}
Note that 
$\pt f$ is 
maximally monotone and
$\sigma$-strongly monotone
by e.g.,
\cite[Theorem~A]{Rock1970},
and
 \cite[Example~22.4(iv)]{BC2017}.
Moreover, we have 
 $(\grad {g})^{-1}=\pt {g}^*$ by
\cite[Remark~on~page~216]{Rock1970}.
%or
%\cite[Th\'{e}or\`{e}me~3.1]{Gossez}.
Therefore, in view of 
\cite[Corollaire~10]{BH77} (see also
\cite{BC2010}) we learn that
\begin{equation}
\text{\label{lem:aux:pddr:b:ii}
${\bf B}$
is maximally monotone and 
${\mu}$-strongly monotone 
and ${\mu}=\min\left\{\sigma,1/\beta\right\}$.
}
\end{equation}

The inclusion in \cref{eq:mono:inc:bAbB} arises
in primal-dual optimality conditions of 
the primal problem \cref{eq:probm} 
and its Fenchel--Rockafellar dual
\cref{eq:probm:FR:dual} given by:
\begin{equation}
 \label{eq:probm}
 \tag{P}
\underset{x\in X}{\operatorname{minimize}}\;\; f(x)+g(Lx)
 \end{equation}

  \begin{equation}
  \tag{D}
 \label{eq:probm:FR:dual}
\underset{y\in Y}{\operatorname{minimize}}\;\; f^*(-L^*y)+g^*(y),
 \end{equation}
 under appropriate assumptions on $f$, $g$
and $L$.

We are now ready for the main result in this section.

\begin{theorem}[\thmtit{application to primal-dual 
Douglas--Rachford method}]
\label{prop:app:DR:pd}
Let 
$\mu= %\gamma^{-1}
\min\left\{\sigma,
1/\beta\right\}$.
Suppose 
$\bf A$ and $\bf B$ are as defined 
in \cref{eq:def:bA}.
Set 
\begin{equation}
T=\tfrac{1}{2}\Big(\Id+R_{  \bf B}
R_{  \bf A}\Big),\quad
 \widetilde{T}=\tfrac{1}{2}\Big(\Id+R_{ \bf A}
 R_{ \bf B}\Big).
 \end{equation}
Let $\bx_0\in X\times Y$, let
$(\bx_n)_\nnn=(T^n \bx_0)_\nnn$,
let
 $(\by_n)_\nnn=(J_{ \bf A}T^n \bx_0)_\nnn$,
 let $(\widetilde{\bx}_n)_\nnn=(\widetilde{T}^n \bx_0)_\nnn$,
and let
 $(\widetilde{\by}_n)_\nnn=(J_{ \bf B}\widetilde{T}^n \bx_0)_\nnn$. 
 Then there exists 
 $\overline{\bx}\in X\times Y
 $,
 $
 \{\overline{\bx}\} =\fix T
 =\fix R_{ \bf B}R_{ \bf A}
$,
there exists 
   $\hat{\bx}\in X\times Y
   $, 
   $\{\hat{\bx}\}=
    \fix R_{ \bf A}R_{ \bf B}=\fix \widetilde{T}$,
 such that $\zer(\bA+\bB)=\{J_{ \bA}\overline{\bx}\}
 =\{J_{ \bB}\hat{\bx}\}$.
 Moreover, the following hold:
 \begin{enumerate}
 \item
 \label{thm:pdDR:i}
 $(\bx_n)_\nnn$ converges strongly to $\overline{\bx}$ with a linear rate 
 $r$, where
\begin{equation} 
\label{eq:rate:r:pd}
%r\in \left]0,
r=
 \tfrac{1}{2(1+\mu)}\Big(\sqrt{2\mu^2+2\mu+1
 +2(1-2(1+\norm{L}^2)^{-1})\mu(1+\mu)}
 +1\Big).
\end{equation} 
 \item
  \label{thm:pdDR:ii}
  $(\by_n)_\nnn$ converges strongly to 
  $J_{ \bA}\overline{\bx}$ with a linear rate 
  $r$ given in \cref{eq:rate:r:pd}.
   \item
 \label{thm:pdDR:iii}
  $(\widetilde{\bx}_n)_\nnn$ converges strongly to 
  $\hat{\bx}$ with a linear rate 
  $r$ given in \cref{eq:rate:r:pd}.
  \item
 \label{thm:pdDR:iv}
   $(\widetilde{\by}_n)_\nnn$ converges strongly to 
  $J_{ \bB}\hat{\bx}$ with a linear rate 
  $r$ given in \cref{eq:rate:r:pd}.
\end{enumerate}
\end{theorem}
{\it Proof}
Note that $\bA$ is 
$\norm{L}$-Lipschitz continuous 
 and $\bB$
is $\mu $-strongly monotone
by 
%\cref{lem:aux:pddr:b}
\cref{lem:aux:pddr:b:i} and \cref{lem:aux:pddr:b:ii}
respectively.
The proof of
\ref{thm:pdDR:i}--\ref{thm:pdDR:iv}
follows from 
 \cref{thm:DR:skew:stmono}
%(respectively \cref{thm:DR:nonex:stmono:linear})
applied with 
$X$ replaced by
$X\times Y$,
$\beta$ replaced by $\norm{L}$,
$A$ replaced by
$ \bA$
 and
$B$ replaced by
$  \bB$.
\qed
\section{Conclusion}
%The main contribution in this paper 
%is the following:
In this paper  we
prove that the Douglas--Rachford method
converges \emph{linearly}
with a \emph{sharp} rate
 when applied to solve \cref{eq:P:sum} in the case $A$ is
Lipschitz continuous 
(but not necessarily a subdifferential operator)
and
$B$ is strongly monotone. 
We also discuss an important application of the results  %\ref{C:R2}
to primal-dual Douglas--Rachford method.
In this case we get sharp rate.
% (see \cref{prop:corres:A:RA}\ref{prop:corres:A:RA:v}).
As a byproduct of our work, we obtain 
useful equivalences between the operator properties
 and the properties of the corresponding resolvent and reflected 
 resolvent.
 
\section*{Acknowledgements}
 This work was done while the authors were visiting 
 the Simons Institute for the Theory of Computing. 
 It was partially supported by the DIMACS/Simons 
Collaboration on Bridging
Continuous and Discrete 
Optimization through NSF grant \# CCF-1740425.

\end{document}